\input amstex

\documentstyle{amsppt}
\document
\topmatter
\title
  Generalized Calabi type K\"ahler surfaces
\endtitle
\author
W\l odzimierz Jelonek, Ewelina Mulawa
\endauthor

\abstract{ In this paper we give a classification  of QCH K\"ahler surfaces of generalized  Calabi type.  }
 \endabstract
\endtopmatter
\define\r{\rightarrow}
\define\G{\Gamma}
\define\DE{\Cal D^{\perp}}

\define\n{\nabla}
\define\om{\omega}
\define\w{\wedge}
\define\k{\diamondsuit}
\define\th{\theta}
\define\p{\partial}
\define\a{\alpha}

\define\lb{\lambda}

\define\1{D_{\lb}}
\define\2{D_{\mu}}
\define\bO{\bar{\Omega}}

\define\0{\Omega}
\define\Om{\Omega}

\define\bJ{\bar J}

\define\J{\Cal J}

\define\De{\Cal D}
\define\R{\Cal R}



{\bf 0. Introduction.}  In  [J-5] the first author has studied QCH semi-symmetric K\"ahler surfaces $(M,g,\bJ)$ which admit an opposite Hermitian structure  which is not locally conformally K\"ahler. A general method of constructing  such manifolds was given there.
 In the present paper we obtain a classification of generalized Calabi type    K\"ahler surfaces.   In particular we obtain, by a different method  than in [A-C-G] a classification of Calabi type K\"ahler surfaces.   We introduce a special orthonormal  frame  $\{E_1,E_2,E_3,E_4\}$  and dual coframe $\{\th_1, \th_2, \th_3, \th_4\}$.   We also find local coordinates $(x,y,z,t)$ in which   $\th_1=fdx$, $\th_2=fdy$  and  $E_4=\p_z, E_3=\beta(z)\p_t$ where $\beta$ is some smooth function.  We show that  such coordinates always exist.   We classify  generalized Calabi type  K\"ahler QCH surfaces which admit an opposite Hermitian structure which is not locally conformally K\"ahler.   The  conditions  on  a Riemannian metric on 4-dimensional manifold   to be   locally conformally K\"ahler  were  studied   for example  in  [D-T].   In our paper we exhibit  a large  class  of   Hermitian  surfaces $(M,J)$ with   degenerate  Weyl  tensor  $W^+$  and  $J$-invariant Ricci tensor, which are not  locally  conformally K\"ahler.
\bigskip
{\bf 1.  Hermitian 4-manifolds. } Let $(M,g,J)$ be an {\it almost Hermitian manifold}, i.e., $J$  is an almost complex structure that is  orthogonal with respect to $g$: $g(X,Y)=g(JX,JY)$ for all $X,Y\in\frak X(M)$. We say that $(M,g,J)$ is a {\it Hermitian  manifold} if the almost Hermitian structure $J$ is integrable, which means that the Nijenhuis tensor $N^J(X,Y)$ vanishes.   This is also equivalent to integrability of the complex distribution $T^{(1,0)}M\subset TM\otimes\Bbb C$.

The {\it K\"ahler form}  of $(M,g,J)$  is $\0(X,Y)=g(JX,Y)$.   If  $d\0=0$  the Hermitian manifold is called a {\it K\"ahler  manifold}.  In what follows we shall consider  K\"ahler manifolds of real  dimension $4$ which are called {\it K\"ahler surfaces}.  Such manifolds are always oriented and we choose an orientation in such a way that the K\"ahler form $\0(X,Y)=g(JX,Y) $ is  self-dual (i.e., $\0\in\w^+M$).

We  investigate K\"ahler surfaces admitting an  opposite  Hermitian structure.  A Hermitian 4-manifold $(M,g,J)$ is said to have an {\it opposite Hermitian structure}  if it admits an orthogonal Hermitian structure $\bar J$ with an anti-self-dual K\"ahler form $\bO$.   A Hermitian manifold $(M,g,J)$ is said to have {\it $J$-invariant Ricci tensor} $\rho$ if $\rho(X,Y)=\rho(JX,JY)$ for all $X,Y\in \frak X(M)$  (see [A-G]).

 {\it Calabi type K\"ahler surfaces} are K\"ahler surfaces which admit a Hamiltonian Killing vector field  X such that  the opposite almost Hermitian structure defined by the complex distribution $\De=\text{ span}\{X,JX\}$ is Hermitian and locally conformally K\"ahler (see [A-C-G]).

{\it  QCH K\"ahler surfaces} are K\"ahler surfaces  $(M,g,\bJ)$
admitting a global, $2$-dimen\-sional, $\bJ$-invariant distribution
$\De$ having the following property: The holomorphic curvature
$K(\pi)=R(X,\bJ X,\bJ X,X)$ of any $\bJ$-invariant $2$-plane $\pi\subset
T_xM$, where $X\in \pi$ and $g(X,X)=1$, depends only on the point
$x$ and the number $|X_{\De}|=\sqrt{g(X_{\De},X_{\De})}$, where
$X_{\De}$ is the orthogonal projection of $X$ on $\De$. In this
case  we have
$$R(X,\bJ X,\bJ X,X)=\phi(x,|X_{\De}|)$$ where $\phi(x,t)=a(x)+b(x)t^2+c(x)t^4$ and
 $a,b,c$ are smooth functions on $M$. Also $R=a\Pi+b\Phi+c\Psi$
 for certain curvature tensors $\Pi,\Phi,\Psi\in \bigotimes^4\frak X^*(M)$
  of K\"ahler type. Every  QCH K\"ahler surface admits an opposite almost Hermitian structure $J$  such that the Ricci tensor $\rho$ of $(M,g,\bJ)$ is $J$-invariant. In  [J-1] it is proved that  every QCH K\"ahler  surface $(M,g,\bJ)$ with  opposite almost Hermitian structure $J$, which is Hermitian and locally conformally K\"ahler, is of Calabi type or  is orthotoric or is hyperk\"ahler.  For a definition of an orthotoric surface see [A-C-G], [J-1].  If $J$ is determined by a foliation into complex curves then  $(M,g,\bJ)$ is of Calabi type.   For example  Eguchi-Hanson surface is of Calabi type  for the standard K\"ahler structure, however is neither of Calabi nor orthotoric type  for infinitely many  K\"ahler structures from its sphere of K\"ahler structures.   The  orthotoric  and Calabi type K\"ahler surfaces were classified  in  [A-C-G].

  By a {\it generalized Calabi type K\"ahler surface} we mean a QCH K\"ahler surface such that the opposite almost Hermitian structure $J$ is determined by a complex foliation into complex curves and is Hermitian.  It just means that  $\De$ is integrable and $J$ is also integrable.

For a Hermitian surface the covariant derivative of the K\"ahler form $\0$ is locally expressed by
$$\n \0=a\otimes\Phi+\J a\otimes\J\Phi, \tag 1.1$$
with $\J a(X)=-a(JX)$,  $\Phi\in LM$, where $LM$ is the bundle of self-dual forms orthogonal to $\0$ and $\J$ is the complex structure on $LM$. The {\it Lee form} $\th$ of $(M,g,J)$ is defined by the equality $d\0=2\th\w \0$. We have $2\th=-\delta\0\circ J$.  We denote by $\tau$ the scalar curvature of  the Riemannian manifold $(M,g)$ i.e.,  $\tau=\text{ tr}_g\rho$.
  An involutive distribution is called a foliation. A foliation $\De$ is called totally geodesic if  every leaf is a totally geodesic submanifold of $(M,g)$, i.e., $\n_XY\in \G(\De)$ for every   $X,Y\in\G(\De)$.

A foliation $\Cal F$
on a Riemannian manifold $(M,g)$ is called {\it conformal} if for every
$V\in\G(T\Cal F)$,
$$L_Vg=\a(V)g$$  on $T\Cal F^{\perp}$, where $\a$ is a
1-form vanishing on $T\Cal F^{\perp}$ (see [J-3]).
\bigskip
{\bf 2. Hermitian surfaces with Hermitian Ricci tensor and  K\"ahler natural opposite almost Hermitian structure.}  Now we recall some results from [J-5].
\par
\bigskip
{\bf Lemma A. }{\it Let $(M,g,J)$ be a Hermitian 4-manifold. Let us assume that  $|\n J| > 0$ on $M$. Then there exists a global oriented orthonormal basis $\{E_3,E_4\}$ of the nullity foliation $\De$, such that for any local orthonormal oriented basis $\{E_1,E_2\}$ of $\DE$
$$\n \0=\a(\th_1\otimes\Phi+\th_2\otimes\Psi), \tag 2.1$$
where  $\Phi=\th_1\w\th_3-\th_2\w\th_4,\ \Psi=\th_1\w\th_4+\th_2\w\th_3$, $\a=\pm\frac1{2\sqrt{2}}|\n J|$ and $\{\th_1,\th_2,\th_3,\th_4\}$ is a cobasis dual to $\{E_1,E_2,E_3,E_4\}$.  Moreover $\delta\0=-2\a\th_3$, $\th=-\a\th_4$.}

Throughout  the paper we shall assume that $|\n J|>0$, hence $\a\ne0$.  Then $\De$ is a  foliation.   If  the natural   opposite almost Hermitian structure is  Hermitian then    $\De$  is totally geodesic  (see [J-4]).
\bigskip

{\bf Lemma B.  }{\it Let $(M,g,J)$ be a Hermitian surface with $J$-invariant Ricci tensor (i.e., $\R(LM)\subset \w^+M$). Let $\{E_1,E_2,E_3,E_4\}$ be a local orthonormal frame such that (2.1) holds. Then
$$\gather\G^3_{11}=\G^3_{22}= E_3\ln \a,  \tag a\\
\G^3_{44}=\G^4_{21}= -\G^4_{12}=- E_3\ln \a, \tag b\\
\G^3_{21}=-\G^3_{12},\ \G^4_{11}=\G^4_{22},  \tag c\\
-\G^3_{21}+\G^4_{22}=\a,\tag d\\
\G^4_{33}=-E_4\ln\a+\a, \tag e\endgather$$
where $\n_XE_i=\sum\om^j_i(X)E_j$ and $\G^i_{kj}=\om^i_j(E_k)$.}

\par
\bigskip
{\bf Lemma C. }{\it Let $(M,g,J)$ be a Hermitian surface   with  J-invariant Ricci tensor. Then $\G^4_{13}=-E_2\ln\a, \ \ \G^4_{23}=E_1\ln\a$ and $d\th$ is anti-self-dual. }

\medskip
If $(M,g,J)$ is a Hermitian surface with $|\n J|>0$ on $M$, then the distributions $\De,\DE$ define a natural opposite almost Hermitian structure $\bar{J}$ on $M$. This structure is defined
as follows $\bar J|_{\De}=-J|_{\De},\bar J|_{\DE}=J|_{\DE}$. In the special basis we just have: $\bar J E_1=E_2, \bar J  E_3=-E_4$.

\par
\bigskip

{\bf Lemma D. } {\it Let $(M,g,J)$ be a Hermitian  4-manifold with Hermitian Ricci tensor  and K\"ahler natural opposite structure. Assume that $|\n J|>0$ on $M$. Then

(a) $\De$ is a totally geodesic foliation,

(b)  $E_3\ln\a=0$,

(c) $\n_{E_4}E_4=0$,

 (d) $d\th=-E_2\a\bar{\Phi}-E_1\a\bar{\Psi}$ where $\bar{\Phi}=\th_1\w\th_3+\th_2\w\th_4,\bar{\Psi}=\th_1\w\th_4-\th_2\w\th_3$.}

\medskip
Hence, if the opposite structure is K\"ahler, then  $E_3\a=0, \G^3_{11}=\G^3_{22}=\G^3_{44}=\G^4_{21}=\G^4_{12}=0$ and
$-\G^3_{21}=\G^3_{12}=\G^4_{11}=\G^4_{22}=\frac12\a$.   We also have  $\G^4_{13}=-E_2\ln\a,\G^4_{23}=E_1\ln\a$.

\bigskip

{\bf Lemma  E.} {\it Let $(M,g,J)$ be a Hermitian  4-manifold with Hermitian Ricci tensor  and K\"ahler natural opposite structure. Assume that $|\n J|>0$ on $M$. If $(M,g,\bJ)$ is a semi-symmetric surface foliated by Euclidean space then  $E_4\ln\a=\frac12\a$.  On the other hand, if  $(Mg,\bJ)$ is a QCH surface with Hermitian opposite structure $J$ such that $\bJ$ is a natural opposite structure for $J$  and $E_4\ln\a=\frac12\a$, then   $(Mg,\bJ)$ is semi-symmetric.}

\medskip
Assume that $\th_1=fdx,\th_2=fdy$. We will later show  that  we can always assume that  $\th_1=fdx,\th_2=fdy$  where  $(x,y,z,t)$  is a local foliated coordinate system, i.e.  $\De=\ker dx\cap \ker dy$.  Then  $\G^1_{32}=\frac\a2,\G^2_{41}=0$.  Then we have  (see [J-5])
$$\gather [E_1,E_4]=-\frac\a2 E_1+E_2\ln\a E_3,\tag2.2\\
[E_2,E_4]=-\frac\a2E_2-E_1\ln\a E_3,\\
[E_1,E_3]=-E_2\ln\a E_4,\\
[E_2,E_3]=E_1\ln\a E_4,\\
[E_3,E_4]=-(-E_4\ln\a+\a)E_3,\\
[E_1,E_2]=\G^1_{12}E_1-\G^2_{21}E_2+\a E_3.\endgather$$

We also have
$$\gather  d\th_1=\G^2_{11}\th_1\w\th_2+\frac12\a\th_1\w\th_4,d\th_2=-\G^1_{22}\th_1\w\th_2+\frac12\a\th_2\w\th_4,\tag2.3\\
d\th_3=-\a\th_1\w\th_2-E_2\ln\a\th_1\w\th_4+E_1\ln\a\th_2\w\th_4+(-E_4\ln\a+\a)\th_3\w\th_4,\\
d\th_4=E_2\ln\a\th_1\w\th_3-E_1\ln\a\th_2\w\th_3.\endgather$$

\medskip
 {\bf Lemma F.} {\it  The two almost Hermitian structures $J,\bJ$  on $M$  such that the $(0,1)$  distribution of $J$ is given by  $\th_1+i\th_2=0,\th_3+i\th_4=0$ and the $(0,1)$ distribution of $\bJ$  by $\th_1+i\th_2=0,\th_3-i\th_4=0$, are integrable  if  equations (2.3) are satisfied.  The form  $d\th=-d(\a\th_4)$ is self-dual with respect to the orientation given by $\bJ$.}

\bigskip
{\bf 3. Special coordinates  for generalized Calabi type   K\"ahler surfaces.}  Let  $(M,g,\bar J)$ be a K\"ahler surface with foliation $\De$  whose leaves are complex curves  such that the natural opposite almost Hermitian structure $J$ determined by $\De$ is Hermitian.   Then $\De$  is conformal and  $d\Omega=2\th\w \Omega$ and  $L_Vg=\th(V)g$ on $\De^{\perp}$ for any $V\in\De$.  It is shown in [J-3] that  $d\th(X,Y)=0$ for all $X,Y\in\De$.

 Note that   Calabi type surfaces $(M,g,J)$ are K\"ahler QCH surfaces with a foliation $\De$ which determines an opposite Hermitian locally conformally K\"ahler structure. Hence $\De$ is conformal and homothetic and we have  $d\th=0$  and hence  $E_1\a=E_2\a=E_3\a=0$.  Hence if $\a$ is not constant  we may assume that  $E_4\a\ne0$   on  $M$: we just  consider the open  submanifold $U$ given by  $U=\{x\in M:E_4\a\ne0\}$.  We construct special coordinates on $M$.
\medskip
{\bf Theorem A.}  {\it Let  $X\in\De$  where  $\De$  is a 2-dimensional foliation on the manifold  $M$, $\dim M=4$. Then there exist foliated coordinates  $(x,y,z,t)$ on $M$  such that  $X=\frac{\p}{\p z}$.}

{\it Proof.}  Assume    $(x',y',z,t)$  are coordinates on $M$ such that  $X=\frac{\p}{\p z}$ and $(dz,dt):\De\r \Bbb R^2$ is a surjection.     As  new coordinates let us take
$(x,y,z,t)$,   where  $(x,y,z',t')$   are  foliated  coordinates, i.e.   $\De=\ker dx\cap \ker dy$.   Then $dx(X)=dy(X)=dt(X)=0, dz(X)=1$  and the 1-forms $dx,dy,dz,dt$ are linearly independent.    Hence $(x, y,$ $z, t)$ is a new foliated coordinate system and in these coordinates  $X=\frac{\p}{\p z}$.$\k$
\medskip
{\bf Theorem B.}  {\it Let  $X,Y\in\De$ be linearly independent vector fields, where  $\De$ is a 2-dimensional foliation on  $M$, dim $M=4$.  Then there exist foliated coordinates  $(x, y, z, t)$ on $M$  such that  $X=f\frac{\p}{\p z},Y=g\frac{\p}{\p t}$.}

{\it Proof.}  Let   $(x',y',z',t)$  be coordinates such that $X=\frac{\p}{\p z'}$ and $(dz',dt):\De\r \Bbb R^2$ is a surjection.    Let   $(x'',y'',z,t')$  be coordinates  such that $Y=\frac{\p}{\p t'}$ and $(dz,dt'):\De\r \Bbb R^2$ is a surjection.   As  new coordinates  take
$(x,y,z,t)$,   where  $(x,y,z_1',t_1')$   are foliated coordinates such that  $\De=\ker dx\cap \ker dy$.   Then $dx(X)=dy(X)=dt(X)=0$  and $dx(Y)=dy(Y)=dz(Y)=0$. The 1-forms $dx,dy,dz,dt$ are linearly independent.     Hence $(x,y,z,t)$ is a new foliated coordinate system and in  the new coordinates    $X=f\frac{\p}{\p z},Y=g\frac{\p}{\p t}$.$\k$

\medskip
{\bf Theorem C.}
{\it Let   $\a\in C^{\infty}(M)$    be a function such that  $X\a\ne0$, where $X\in\De$.  Then there exist foliated coordinates  $(x,y,z,t)$  for which   $\a=\a(z)$ and $X=f\frac{\p}{\p z}$ for a certain function $f$.}

{\it Proof.}   Take coordinates $(x,y,z,t)$ given by Th. A and let  $x'=x$, $y'=y$, $z'=g(\a(x,y,z,t))$, $t'=t$  for some smooth function $g$ on $\Bbb R$ with $g'\ne0$.  Then $\a$ is a function of $z'$ only and  $X=\frac{\p}{\p z}= g'(\a(x,y,z,t))\frac{\p \a}{\p z}(x,y,z,t)\frac{\p}{\p z'}= f\frac{\p}{\p z'}$  where  $$f=g'(\a(x,y,z,t))\frac{\p \a}{\p z}(x,y,z,t).$$ $\k$

\medskip
{\bf Theorem  D.}  {\it Let    $a,b:\Bbb R^4\r \Bbb R$  be smooth functions such that  $\p_ya=\p_xb$.  Then there exists a smooth function  $f:\Bbb R^4\r \Bbb R$ such that  $\p_xf=a,\p_yf=b$.}

  {\it Proof.} Define   $f(x,y,z,t)=A(z,t)+\int_{\gamma}a(x,y,z,t)dx+b(x,y,z,t)dy$   where $A$ is some smooth function and  $\gamma$  is a smooth curve contained in the plane $z=\text{const},t=\text{const}$ which joins a point  $(0,0,z,t)$  to a point  $(x,y,z,t)$.  Then   $f$  is well defined since the form   $\om=adx+bdy$  is closed (hence exact)  in the plane  $z=\text{const},t=\text{const}$.  It is easy to see that $\p_xf=a,\p_yf=b$.$\k$
\bigskip

{\bf Lemma  G.}  {\it  Let  $\th$ be a  1-form on $\Bbb R^4$ such that $d\th(\p_x,\p_y)=0$.  Then there exists a function $f:\Bbb R^4\r \Bbb R$ such that  $\th(\p_x)=\p_xf,\th(\p_y)=\p_yf$.}

{\it Proof.} Let   $a=\th(\p_x),b=\th(\p_y)$.   Then  $\p_ya=\p_xb$ and we can apply Th. D.$\k$
\medskip
{\bf Theorem E.} {\it Let  $(M,g,J)$  be a K\"ahler surface with conformal 2-dimensional foliation $\Cal F$  such that  $d\th(X,Y)=0$ for  $X,Y\in T\Cal F$.
Then for horizontal vector fields $Z,T\in (T\Cal F)^{\perp}$ we have   $g(T,Z)=fp^*h(T,Z)$   where $p:M\r N$  is a local submersion, whose fibers are leaves of the foliation   $\Cal F$ and $f$  is a certain function on  $M$.}

  {\it Proof.}  Let  $L_Vg(Y,Z)=\th(V)g(Y,Z)$.  We have   $\th(X)=d\ln f(X)$  for all $X\in T\Cal F$ for a certain positive function  $f$ (see Lemma G).   Define a metric  $h$  on $N$ by the formula  $h(X, Y)=\frac1fg(X^*,Y^*)$  where $X^*,Y^*$  are horizontal  lifts of  $X,Y$, i.e.   $p(X^*)=X$,  $p(Y^*)=Y$  and $X^*,Y^*\in (T\Cal F)^{\perp}$.  We will show that $h$ is well defined.   It is enough to show that for $V\in T\Cal F$   we have  $V(\frac1fg(X^*,Y^*))=0$.  But   $$\gather V\big(\frac1fg(X^*,Y^*)\big)= \\ -\frac{Vf}{f^2}g(X^*,Y^*)+\frac1f(L_Vg)(X^*,Y^*)\\+\frac1fg([V,X^*],Y^*)+\frac 1fg(X^*,[V,Y^*])\\=-\frac{(d\ln f) V}fg(X^*,Y^*)+\frac1f (d\ln f) V g(X^*,Y^*)=0\endgather$$  since the fields $[V,X^*],[V,Y^*]$ are tangent to the leaves of the foliation.$\k$

\medskip
Remark.   In the case of QCH  K\"ahler surfaces  with Hermitian opposite  almost Hermitian structure given by a complex foliation $\Cal F$ we always have  $d\th(X,Y)=0$ for  $X,Y\in T\Cal F$  (see [J-3], Remark 3.3,  p.  234).

\medskip
{\bf Theorem F.}  {\it Let  $(M,g,J)$  be a QCH K\"ahler surface with a complex foliation $\De$  determining a Hermitian opposite natural structure.   Then  there exists an orthonormal frame $\{E_1,E_2\}$  on $\De^{\perp}$ such that $\th_1=fdx,\th_2=fdy$ where $f$ is a positive function and $(x,y,z,t)$ is a foliated coordinate system on $M$.}

{\it Proof.}   Let us take  isothermal coordinates on $N$  for $h$.   Then   $h=k(dx^2+dy^2)$   and   $g(Y,Z)=f_1k(dx^2+dy^2)$ for  $Y,Z\in T\De^{\perp}$ (see Th. E).    Now take  $f=\sqrt{f_1k}$$\k$
\medskip
Now we show  that  for Calabi surfaces  with non constant $\alpha$  we can   find local coordinates for which     $E_1=\frac1f\p_x+l\p_t,E_2=\frac1f\p_y+n\p_t,E_3=\beta\p_t,E_4=\p_z$  where $g=\beta g_1$ $\a\beta=e^A$,  $A=\int \a$  and  $g_1=g_1(x,y,t)$.  We assume that  $E_4\a\ne 0$  on $M$.   Note that in our case  $E_1\a=E_2\a=E_3\a=0$,  hence we assume that $\a$ is nonconstant.  From the above theorems it follows that  we can take foliated coordinates such that  $E_4=h\p_z$ and $E_3\in \text{span}(\p_z,\p_t)=\De$.   Since $E_3\a=0$  and $\a_z\ne0$  it follows that  $E_3=g\p_t$ for some function $g(x,y,z,t)$.

 Now we show that we can take  $h=1$.    Note that $\ker dz=\text{ span}(E_1,E_2,E_3)=\text{span} (\p_x,\p_y,\p_t)$.  $E_4\a$  depends only on $z$  since  using  the equations for Lie brackets we obtain  $E_1E_4\a=[E_1,E_4]\a=-\frac{\a}2E_1\a=0$ and so on.   This implies  $E_4\a=\phi(z)$.   Now   $E_4\a=h\a'(z)=\phi(z)$.  Hence   $h$ is a function of $z$ only   and making an appropriate transformation of coordinates we can assume that   $h=1$.

 Now  from  $[E_3,E_4]=-(-\frac{\a'}{\a}+\a)E_3$  it is easy to check that  $g=\beta g_1(x,y,t)$  where  $\a\beta=e^A$.    Now  we change coordinates  by $x'=x, y'=y, t'=\int\frac1{g_1(x,y,t)}dt, z'=z$.  Then   $\p_{z'}=\p_z, \p_{t'}=g_1\p_t$.  Then in the new coordinates   $g_1=1$  and  $E_4=\p_{z'},E_3=\beta(z')\p_{t'}$. We write these new coordinates as $(x,y,z,t)$. Note that this system of coordinates is foliated.   In the case  $\a=\text{ const}\ne0$   we use Th. B  to find foliated coordinates such that  $E_4=h\p_z,E_3=g\p_t$.   Now it is easy to see that $\p_th=0$.   Take new coordinates $x'=x,y'=y,t'=t,z'=\int \frac1{h(x,y,z)}dz$.    Then in the new coordinates we have  $E_4=\p_z$, $E_3=g\p_t$.  Now again  $g=\beta(z)g_1(x,y,t)$.    Using the transformation $x'=x, y'=y, t'=\int\frac1{g_1(x,y,t)}dt, z'=z$ we can assume  that  $E_4=\p_z$, $E_3=\beta(z)\p_t$  and  the new coordinates are foliated.  After some calculations one can prove that if $\alpha=\text{ const}$  we can always find coordinates in which  $E_1=\frac1f\p_x+l\p_t, E_2=\frac1f\p_y+n\p_t, E_3=\beta\p_t, E_4=\p_z$, so this case is the case of Calabi surfaces.  In fact  then  $d\th_4=0$  and hence  $\th_4 =dz+dF$  where  $F=F(x,y)$.   Taking transformation  $x'=x, y'=y, z'=z+F(x,y),t'=t$  we have foliated coordinates and in these coordinates  $k,m=0$  where  $E_1=\frac1f\p_x+k\p_z+l\p_t,E_2=\frac1f\p_y+m\p_z+n\p_t,E_3=\beta\p_t,E_4=\p_z$.

Now we consider the case of QCH  K\"ahler surfaces  with  opposite Hermitian structure which is not locally conformally K\"ahler.   Take foliated coordinates  such that $\th_1=fdx$, $\th_2=fdy$ and  $E_1=\frac1f\p_x+k\p_z+l\p_t,E_2=\frac1f\p_y+m\p_z+n\p_t, E_4=r\p_z,  E_3=g\p_t$.  It is easy to see that  $\p_tr=0$  hence   $r=r(x,y,z)$.    Let us change coordinates by $x'=x,y'=y,z'=\int\frac1{r(x,y,z)}dz,t'=t$.
Then in the new coordinates   $\a=\a(x,y,z)$  and  $E_4=\p_z,E_3=g(x,y,z,t)\p_t$.  Let  $\beta(x,y,z)$  satisfy  $\p_z\ln\beta+\p_z\ln\alpha=\alpha$. Write    $A=A(x,y,z)=\int \a dz$.  Hence  $\a\beta=e^A$.  Note that   $\p_z\ln g=\p_z\ln\beta$ and hence  $g=\beta g_1(x,y,t)$.    Now changing coordinates by $x'=x, y'=y, z'=z, t'=\int\frac1{g_1(x,y,t)}dt$  we obtain new foliated coordinates in which   $E_4=\p_z,E_3=\beta\p_t$.  Since  $\th_1=fdx$, $\th_2=fdy$, from our equations we get

$$\big[\frac 1f\p_x+k\p_z+l\p_t,\p_z\big]=-\frac{\alpha}2\big(\frac 1f\p_x+k\p_z+l\p_t\big)+\big(\frac1f\p_y\ln\a+m\p_z\ln\a\big)\beta\p_t.\tag 3.1$$
Hence  $\p_zf=-\frac{\a}2f, \p_zk= \frac{\a}2k$  and $f=e^{-\frac A2}h(x,y)$  and  $k=e^{\frac A2}k_1(x,y,t)$. Analogously $\p_zm=\frac{\alpha}2m$ and  $m=e^{\frac A2}m_1(x,y,t)$. We also have

$$\big[\frac 1f\p_x+k\p_z+l\p_t,\beta\p_t\big]=-\big(\frac 1f\p_y\ln\a+m\p_z\ln\a\big)\p_z.$$
 This implies    $\p_tk=\frac1{f\beta}\p_y\ln\a+m\frac{\p_z\ln\a}{\beta}$.  Consequently, we obtain  $$\p_{zt}k= \frac{\alpha}2k_t =\frac{\alpha}2(\frac1{f\beta}\p_y\ln\a+m\frac{\p_z\ln\a}{\beta}).\tag 3.2$$  On the other hand,  $$\p_{zt}k=\big(\frac1{f\beta}\p_y\ln\a\big)_z+m\a\frac{\p_z\ln\a}{2\beta}+m\big(\frac{\p_z\ln\a}{\beta}\big)_z.\tag 3.3$$
Note also that  $\p_tk=\frac1{f\beta}\p_y\ln\a+m\frac{\p_z\ln\a}\beta$ and  $\p_tm=-\frac1{f\beta}\p_x\ln\a-k\frac{\p_z\ln\a}\beta$.  Hence if $m_t=0$  then $k_t=0$  and  $m=-\frac1{f\p_z\ln \alpha}\p_y\ln\alpha,k=-\frac1{f\p_z\ln \alpha}\p_x\ln\a$.  It is easy to check  that in this case  $E_1\ln\a=E_2\ln\a=0$  and consequently $d\th=0$,  which means that the opposite structure is locally conformally K\"ahler.

Hence we can assume $m_t\ne0,k_t\ne0$.
Since $f,\beta,\a$ do not depend on $t$ and $m$ is a nonconstant function depending on $t$ we get  $(\frac{\p_z\ln\a}{\beta})_z=0$.  It  follows that  $\beta= C(x,y)(\ln\a)_z$, where $C\ne0$.   Thus  since  $\p_z(\ln\beta)=-\p_z\ln\a+\a$  we obtain $\a''=\a\a'$ and   $\a'=\frac12\a^2+D(x,y)$.   Now in the case  $D>0$  we get   $\a=2a\tan(a (z +\phi(x,y)) )$  where  $a=a(x,y)=\sqrt{\frac D2}$. Taking new coordinates with  $z'=z+\phi$  we can assume that   $\a=2a\tan(a z  )$.  In the case   $D<0$  we get
$\a=-2a\coth(a z)$ or $\a=-2a\tanh(az)$, where $a=a(x,y)=\sqrt{-\frac D2}$.  If   $D=0$  in some open subset of $M$  we get  $\a=-\frac2{z+\phi(x,y)}$  and changing  coordinates  by  $x'=x,y'=y,z'=z+\phi(x,y),t'=t$  we can assume that   $\a=-\frac2z$.  From (3.2) and (3.3) it follows that   $$\frac1{f\beta}\p_y\ln\a=e^{\frac A2} c(x,y)\tag 3.4$$ for some function $c$.   We shall consider the case  $\a=2a\tan(a z  )$ . Then   $$\a_y= 2a_y\tan az+2a\frac1{\cos^2az}a_yz\tag 3.5$$  and   $(\ln\a)_y= (\ln a)_y+\frac{a_yz}{\sin az\cos az}$.  Note that $\beta=C(x,y)\frac{2a}{\sin2az}$ and  we get $$\frac {\p_y\ln\a}{h\beta}=\frac{a_y}{Ch}\big(\frac1{2a^2}\sin2az+\frac{ z}{a}\big)=c(x,y).\tag 3.6$$  From  (3.6) we obtain   $a_y=0$.  Similarly we show that $a_x=0$  and hence  $a=\text{ const}$.    Now change coordinates as follows: $x'=x,y'=y,z'=z,t'=\frac1{C(x,y)}t$.  These are foliated coordinates and  $\p_{t'}=C(x,y)\p_t, \p_{z'}=\p_z$.  This yields   $\beta=\frac{2a}{\sin2az}$.  The case $D<0$ is similar. This implies that $\alpha,\beta$ are functions of $z$ only in the introduced coordinates and the only possibilities for the function $\alpha$ are  $\a=-\frac2z,\a=2a\tan az,\a=-2a\coth az, \a=-2a\tanh az$, where $a\in\Bbb R,a\ne0$.

\medskip
{\bf 4.  Calabi type K\"ahler surfaces.}
 Next   using our method  we classify  Calabi type K\"ahler surfaces  already classified in [A-C-G].  Let  $\a(z)$  be any smooth, non-vanishing function defined on an open subset $V\subset \Bbb R$ and $A=\int \a,  \beta \alpha=e^A$.
\medskip
{\bf Theorem 1.}  {\it  Let  $U\subset \Bbb R^2$ be an open set and let $g_{\Sigma}=h^2(dx^2+dy^2)$ be a Riemannian metric on $U$, where $h:U\r \Bbb R$ is a positive function $h=h(x,y)$.  Let $\om_{\Sigma}=h^2dx\w dy$ be the volume form of $\Sigma=(U,g_{\Sigma})$. Let $M=U\times N$, where $N=V\times\Bbb R$. Let us define a metric $g$  on $M$  by $g(X,Y)=e^{-A}g_{\Sigma}(X,Y)+\th_3(X)\th_3(Y)+\th_4(X)\th_4(Y)$,   where $\th_3=\frac 1{\beta}(dt-l_2(x,y)dx-n_2(x,y)dy)$,  $A=\int \a$,  $\beta \alpha=e^A$ and
$\th_4=dz$ and $d(l_2dx+n_2dy)=\om_{\Sigma}$.  Then   $(M,g)$ admits a K\"ahler structure $\bJ$ with the K\"ahler form $\bO=e^{-A}\om_{\Sigma}+\th_4\w \th_3$  and a Hermitian structure $J$ with the K\"ahler form $\Om=e^{-A}\om_{\Sigma}+\th_3\w \th_4$.  The Ricci tensor of $(M,g)$ is $J$-invariant and $J$ is  locally conformally K\"ahler.  The Lee form of $(M,g,J)$ is $\th=-\a\th_4$.  The scalar curvature of $(M,g)$ is $\tau=2(-\frac{(\Delta\ln h)e^A}{h^2}-2\a^2+2\a'+\frac{\beta''}{\beta}-2(\frac{\beta'}{\beta})^2)$.}
\medskip
{\it Proof.}   We shall use equations (2.2).
Let us take a coordinate system such that  $\th_1=fdx,\th_2=fdy$,  $\a=\a(z)$ and assume  $E_1\a=E_2\a=0$.  Then  if  $\alpha$ is not constant we obtain  $E_1=\frac1f\p_x+l\p_t, E_2=\frac1f\p_y+n\p_t, E_3=\beta\p_t, E_4=\p_z$. However if $\alpha$ is constant we can change coordinates such that in new coordinates we have also $E_1=\frac1f\p_x+l\p_t, E_2=\frac1f\p_y+n\p_t, E_3=\beta\p_t, E_4=\p_z$.  This implies
 $\th_1=fdx$, $\th_2=fdy$, $\th_4=dz$, $\th_3=\frac1{\beta}( dt- lfdx-nfdy)$. Then $E_4\ln\a=\frac{\a'}{\a}$ and
$$E_1\ln\a=0,   \   \   E_2\ln\a=0.\tag4.1$$
We have $[E_1,E_4]=\frac{f_z}{f^2}\p_x-l_z\p_t=-\frac\a{2f}\p_x-\frac\a2 l\p_t$. Hence $
l_z=\frac{\a}2 l, f_z=-\frac{\a}2 f$. This implies   $f=e^{-\frac A2}h(x,y)$.  On the other hand, $[E_1,E_3]=-\beta l_t\p_t=0$  and
$$l_t=0.\tag 4.2$$
This yields
$$l=l_1(x,y)e^{\frac A2}.\tag4.3$$
Similarly $[E_2,E_3]=-gn_t\p_t=0$
 and consequently, $n_t=0$.

Since $[E_2,E_4]=-n_z\p_t+\frac{f_z}{f^2}\p_y=-\frac\a{2f}\p_y-\frac12\a n\p_t$, we get
$n_z=\frac12\a n$. Hence
$$n=n_1(x,y)e^{\frac A2}.\tag4.4$$
 Therefore
$\th_3=\frac 1{\beta}(dt-l_2(x,y)dx-n_2(x,y)dy)$
and
$\th_4=dz$.

Now we prove that
$d\th_3=-\a\th_1\w\th_2+(-\frac{\a'}{\a}+\a)\th_3\w\th_4=-\a\th_1\w\th_2+\frac{\beta_z}{\beta}\th_3\w\th_4$ if    $d(l_2dx+n_2dy)=\om_{\Sigma}$.

In fact,
$$\gather d\th_3=d\big(\frac 1{\beta}(dt-l_2(x,y)dx-n_2(x,y)dy)\big)=\\ -\a\th_1\w\th_2+\frac{\beta_z}{\beta}\th_3\w\th_4 \endgather$$
after some easy computation if we assume that   $d(l_2dx+n_2dy)=\om_{\Sigma}=h^2dx\w dy$.

Now we have   $d(\th_1\w\th_2)=d(e^{-A}h^2dx\w dy)=-e^{-A}\a dz\w h^2dx\w dy=-\a dz\w\th_1\w\th_2=-\a\th_4\w\th_1\w\th_2$. From (2.3) it is also clear that  $d(\th_3\w\th_4)=d(\th_3)\w\th_4= -\a\th_1\w\th_2\w\th_4 $.   Hence  $d\bO=0$, where $\bO(X,Y)=g(\bJ X,Y)$, and $\bO=\th_1\w\th_2-\th_3\w\th_4$. Hence in view of Lemma F, $(M,g,\bJ)$ is a K\"ahler surface and the Lee form of $(M,g,J)$ is $\th=-\a\th_4$.

Now we show that $\p_t$ is a real holomorphic vector field.  Note that the opposite K\"ahler structure $\bJ$ satisfies  $\bJ\p_z=\beta\p_t,\bJ\p_t=-\frac1{\beta}\p_z,  \bJ\p_x=\p_y+fn\p_t+\frac1{\beta} fl\p_z, \bJ\p_y=-\p_x-fl\p_t+\frac1{\beta} fn\p_z$.
It is clear that $L_{\p_t}\bJ=0$.  The vector field  $X=\p_t+i\frac1{\beta}\p_z$ is  holomorphic.  The form  $\psi=\frac12(dt-i\beta dz)$  satisfies  $\psi(X)=1$ and  $\psi(\p_t-i\frac1{\beta}\p_z)=0$.   If   $\Psi$ is a holomorphic $(1,0)$ form such that  $\Psi(X)=1$,   then  $\Psi-\psi=0 \text{  mod }\{dx,dy\}$.  If  $\bO$ is the K\"ahler form of the K\"ahler manifold $(M,\bJ,g)$, then
$\frac12\bO^2=e^{-A}h^2\frac1\beta dx\w dy\w dz\w dt$.  We also have  $\frac{i^2}4(dx+idy)\w\overline{(dx+idy)}\w\Psi\w \overline{\Psi}=\frac14{\beta} dx\w dy\w dz\w dt$.   Hence the Ricci form  of  $(M,g,\bJ)$  is   $$\gather \rho=-\frac12dd^c\ln (e^{-A})
-dd^c\ln h+dd^c\ln \beta\\= \frac12dd^cA-d((\ln h)_xdy-(\ln h)_ydx)+  d(\frac{\beta'}{\beta}\th_3)\\=-\Delta\ln h \frac{e^A}{h^2}\th_1\w\th_2+(\a'-\frac32\a^2)\th_1\w\th_2+(\frac{\beta''}{\beta}-2(\frac{\beta'}{\beta})^2
+\frac12\a'-\frac{\a\beta'}{2\beta})\th_4\w\th_3.\endgather$$    Consequently, $E_1,E_2,E_3,E_4$ are eigenfields of the Ricci tensor.
Since
$$\gather \rho=(-(\Delta\ln h)e^A\frac1{h^2}-\frac32\a^2+\a')\th_1\w\th_2\\+(\frac12\a'+\frac{\beta''}{\beta}-2(\frac{\beta'}{\beta})^2
-\frac{\beta'\a}{2\beta})\th_4\w\th_3,\endgather$$ it follows that $\frac12\tau=-\frac{(\Delta\ln h)e^A}{h^2}-2\a^2+2\a'+\frac{\beta''}{\beta}-2(\frac{\beta'}{\beta})^2$.$\k$

\medskip
{\bf Remark. }   It is easy to show that  $X=\p_t$  is a Hamiltonian Killing  vector field.  It  also has a special K\"ahler Ricci potential in the sense of Derdzi\'nski-Maschler. ([D-M-1], [D-M-2]). In fact
$\bar J X=-\frac1{\beta}E_4=-\n u$  where $u$ is a potential of $X$, so $X=\bJ\n u$ and   $u(z)=\int \frac1{\beta}dz$.   Now since $\n_{E_4}E_4=0,\n_{E_3}E_4=-(-\frac{\a'}\a+\a)E_3$  it follows that  $X,JX$  are eigenfields of the Hessian  $H^u$.
\bigskip

{\bf 5.  Generalized Calabi type K\"ahler surfaces.}    The semi-symmetric K\"ahler surfaces of this type are classified in [J-5].   Note that we can take  $kf=H\sin\frac12(t+\phi(x,y)),  mf=H\cos\frac12(t+\phi(x,y))$  for a certain function $\phi$  and changing coordinates $x'=x, y'=y, z'=z, t'=t+\phi(x,y)$  we can assume that
$kf=H\sin\frac12t$, $fm=H\cos\frac12t$.  We also have in these coordinates $\a=-\frac2z$.  Now we  classify the remaining cases.  First we assume  $\a=2a\tan az, a\in\Bbb R,a\ne0$.

{\bf Theorem  2.}  {\it  Let  $U\subset \Bbb R^2$ and let $g_{\Sigma}=h^2(dx^2+dy^2)$ be a Riemannian metric on $U$, where $h:U\r \Bbb R$ is a positive function $h=h(x,y)$.  Let $\om_{\Sigma}=h^2dx\w dy$ be the volume form of $\Sigma=(U,g)$. Let $M=U\times N$, where $N=\{(z,t)\in \Bbb R^2:|z|<\frac{\pi}{2|a|}\}$. Define a metric $g$  on $M$  by $g(X,Y)=(\cos az)^2g_{\Sigma}(X,Y)+\th_3(X)\th_3(Y)+\th_4(X)\th_4(Y)$,   where $ \th_3=\sin 2azdt-(\cos 2at\cos (2az)H(x,y)+\sin (2az)l_2(x,y))dx-(-\sin 2at\cos (2az)H(x,y)+\sin (2az)n_2(x,y))dy$,
$$\th_4=dz-\sin 2atH(x,y)dx-\cos 2atH(x,y)dy$$ and the function $H$ satisfies the equation  $\Delta \ln H=(\ln H)_{xx}+(\ln H)_{yy}=2a^2h^2-4a^2H^2$ on $U$, $l_2=-\frac1{2a}(\ln H)_y, n_2=\frac1{2a}(\ln H)_x$.  Then   $(M,g)$ admits a K\"ahler structure $\bJ$ with the K\"ahler form $\bO=(\cos az)^2\om_{\Sigma}+\th_4\w \th_3$  and a Hermitian structure $J$ with the K\"ahler form $\Om=(\cos az)^2\om_{\Sigma}+\th_3\w \th_4$.  The Ricci tensor of $(M,g)$ is $J$-invariant and $J$ is not locally conformally K\"ahler.  The Lee form of $(M,g,J)$ is $\th=-\a\th_4$,  where $\a=2a\tan az$.  }

\medskip

{\it Proof.} Let us take a coordinate system such that $E_1=\frac1f\p_x+k\p_z+l\p_t,E_2=\frac1f\p_y+m\p_z+n\p_t,E_3=\frac1\beta\p_t,E_4=\p_z$. Then  $\th_1=fdx, \th_2=fdy, \th_4=dz-(fk)dx-(fm)dy,\th_3=\beta dt-(\beta lf)dx-(\beta nf)dy$.  Let  $\a=2a\tan az,\beta=\sin 2az$. Then

$$E_1\ln\a=\frac{2a k}{\sin 2az},\  \  E_2\ln\a=\frac {2am}{\sin 2az}.\tag5.1$$
We have $[E_1,E_4]=\frac{f_z}{f^2}\p_x-k_z\p_z-l_z\p_t=-\frac\a{2f}\p_x-\frac\a2k\p_z-\frac\a2 l\p_t+\frac {2am}{\sin^22az}\p_t$. Hence $k_z=a(\tan az) k,
l_z=a(\tan az) l-\frac{2am}{sin^22az}, f_z=-a\tan azf$. This implies   $$f=\cos azh(x,y), k=\frac{k_1(x,y,t)}{\cos az}.$$ On the other hand, $[E_1,E_3]=-\frac1\beta( k_t\p_z+ l_t\p_t)-k\frac{\beta_z}{\beta^2}\p_t=-\frac{2a m}{\sin 2az}\p_z$  and
$$l_t= -2ak\cot 2az, k_t=2am.\tag 5.2$$
This yields
$$l=m\cot 2az+\frac{l_1(x,y)}{\cos az}.\tag5.3$$
Similarly $[E_2,E_3]=-m\frac{\beta_z}{\beta^2}\p_t-\frac1{\beta} (m_t\p_z+ n_t\p_t)=\frac {2ak}{\sin 2az}\p_z$, what implies
$$m_t=-2ak,n_t= -2am\cot 2az.$$
Since $[E_2,E_4]=-m_z\p_z-n_z\p_t+\frac{f_z}{f^2}\p_y=-\frac\a{2f}\p_y-\frac12\a m\p_z-\frac12\a n\p_t-\frac {2a k}{\sin^22az}\p_t$, we get
$m_z=a\tan az m$.  Hence  $m=\frac{m_1(x,y,t)}{\cos az}$ and
$$n=-k\cot 2az+\frac{n_1(x,y)}{\cos az}.\tag5.4$$
Let us take  coordinates in which  $kf=\sin 2atH(x,y), mf=\cos 2at H(x,y)$,  where  $f=\cos azh(x,y)$ and $\th_1\w\th_2=(\cos az)^2h^2dx\w dy$.  Then
$$lf=\cot 2az\cos 2atH(x,y)+l_2(x,y),nf=-\cot 2az\sin 2atH(x,y)+n_2(x,y).$$  Hence
$$\gather \th_3=\sin 2azdt-(\cos 2at\cos 2azH(x,y)+\sin 2azl_2(x,y))dx-\\(-\sin 2at\cos 2azH(x,y)+\sin 2azn_2(x,y))dy\endgather$$
and
$\th_4=dz-\sin 2atH(x,y)dx-\cos 2atH(x,y)dy$.

Now we prove that
$d\th_3=-\a\th_1\w\th_2-E_2\ln\a\th_1\w\th_4+E_1\ln\a\th_2\w\th_4+(-E_4\ln\a+\a)\th_3\w\th_4$ and
$d\th_4=E_2\ln\a\th_1\w\th_3-E_1\ln\a\th_2\w\th_3$  if $l_2=-\frac1{2a}(\ln H)_y, n_2=\frac1{2a}(\ln H)_x$, where  $\Delta\ln H=-4a^2H^2+2a^2h^2$.

In fact,  $$\gather-\a\th_1\w\th_2-E_2\ln\a\th_1\w\th_4+E_1\ln\a\th_2\w\th_4+(-E_4\ln\a+\a)\th_3\w\th_4\\=-2a\tan az f^2dx\w dy -\frac{2amf}{\sin 2az}dx\w(dz-\sin 2atHdx-\cos 2atHdy)\\+\frac{2akf}{\sin 2az}dy\w(dz-\sin 2atHdx-\cos 2atHdy)\\- \frac{2a\cos 2az}{\sin 2az}(\sin 2azdt-\cos 2az \cos 2atH dx\\-\sin 2az l_2 dx+\sin  2at\cos 2az H dy -\sin 2az n_2dy)\\ \w(dz-\sin 2atHdx-\cos 2atHdy)=-a\sin 2azh^2dx\w dy\\+ \frac{2a}{\sin 2az}\cos^22atH^2dx\w dy -\frac{2a}{\sin 2az}\sin^22atH^2dy\w dx +\frac{2a\sin 2atH}{\sin 2az}dy\w dz\\-\frac{2a\cos 2atH}{\sin 2az}dx\w dz-\frac {2a\cos 2az}{\sin 2az}(\sin 2azdt\w dz\\-\sin 2az\sin 2atHdt\w dx-\sin 2az\cos 2atHdt\w dy\\-\cos 2az\cos 2atHdx\w dz-\sin 2azl_2dx\w dz+\cos 2az\cos^22atH^2dx\w dy\\ +\sin 2azl_2\cos 2atHdx\w dy+\cos 2az\sin 2atHdy\w dz\\-\sin 2azn_2dy\w dz-\cos 2az\sin^22atH^2dy\w dx +n_2\sin 2az\sin 2atHdy\w dx). \endgather$$

On the other hand,
$$\gather d\th_3= 2a\cos 2azdz\w dt +2a\sin 2az\cos 2atHdz\w dx+2a\cos 2az\sin 2atHdt\w dx\\ -\cos 2az\cos 2atH_ydy\w dx-2al_2\cos 2azdz\w dx-\sin 2azl_{2y}dy\w dx\\ +2a\cos 2az\cos 2at Hdt\w dy-2a\sin 2at\sin 2azHdz\w dy\\+\sin 2at\cos 2azH_xdx\w dy-2a\cos 2azn_2dz\w dy-\sin 2azn_{2x}dx\w dy.\endgather$$

It is clear that $\th_3$ satisfies (2.3) if $l_2=-\frac1{2a}(\ln H)_y, n_2=\frac1{2a}(\ln H)_x$ and  $\Delta\ln H=2a^2h^2-4a^2H^2$, where $\Delta f=f_{xx}+f_{yy}$.

Now  $$\gather E_2\ln\a\th_1\w\th_3-E_1\ln\a\th_2\w \th_3\\=\frac{2a\cos 2atH}{\sin 2az}dx\w(\sin 2az dt+\sin 2at\cos 2az Hdy- \sin 2azn_2dy)\\-\frac{2a\sin 2at H}{\sin 2az} dy \w(\sin 2azdt-\cos 2az\cos 2at Hdx-\sin 2azl_2dx)\\=2a\big(\cos 2atHdx\w dt-\cos 2atHn_2dx\w dy +\frac{\cos 2at\sin 2at\cos 2azH^2}{\sin 2az}dx\w dy\\-\sin 2at Hdy\w dt+\frac{\sin2a t \cos 2at\cos 2azH^2}{\sin 2az}dy\w dx+\sin 2atHl_2dy\w dx\big).\endgather$$

On the other hand,

$d\th_4=-2a\cos 2atHdt\w dx-\sin 2atH_ydy\w dx+2a\sin 2atHdt\w dy-\cos 2atH_xdx\w dy$.

It is clear that $\th_4$ satisfies (2.3) if $l_2=-\frac1{2a}(\ln H)_y, n_2=\frac1{2a}(\ln H)_x$.

Now we have   $d(\th_1\w\th_2)=d((\cos az)^2h^2dx\w dy)=-2a\sin az\cos azdz\w h^2dx\w dy=-2a\tan azdz\w\th_1\w\th_2=-\a\th_4\w\th_1\w\th_2$. From (2.3) it is also clear that  $d(\th_3\w\th_4)=-\a\th_4\w\th_1\w\th_2$.   Hence  $d\bO=0$, where $\bO(X,Y)=g(\bJ X,Y)$ and $\bO=\th_1\w\th_2-\th_3\w\th_4$. Hence in view of Lemma F $(M,g,\bJ)$ is a K\"ahler surface and the Lee form of $(M,g,J)$ is $\th=-\a\th_4$.

Now we show that $\De=\text{ span}\{E_3,E_4\}$ is a conformal foliation.  Note that the opposite K\"ahler structure $\bJ$ satisfies  $\bJ\p_z=\frac1\beta\p_t,\bJ\p_t=-\beta\p_z,  \bJ\p_x=\p_y+fm\p_z+fn\p_t-\frac1\beta fk\p_t+\beta fl\p_z, \bJ\p_y=-\p_x-fk\p_z-fl\p_t-\frac1\beta fm\p_t+\beta fn\p_z$.  Hence  $(L_{\p_t}\bJ)\p_x,(L_{\p_z}\bJ)\p_x,(L_{\p_t}\bJ)\p_y,(L_{\p_z}\bJ)\p_y\in\De$,  which means that $L_{\xi}\bJ(\DE)\subset \De$ for $\xi\in\G(\De)$. Thus the foliation $\De$ is conformal and since  $d\th^-=0$, it follows from [J-3] that $(M,g,\bJ)$ is a QCH K\"ahler surface.$\k$
\bigskip
Next we consider the remaining cases  $\a=-2a\coth az, a\in\Bbb R,a\ne0$  and  $ \a=-2a\tanh az, a\in\Bbb R,a\ne0$.
\bigskip
{\bf Theorem  3.}  {\it  Let  $U\subset \Bbb R^2$ and let $g_{\Sigma}=h^2(dx^2+dy^2)$ be a Riemannian metric on $U$, where $h:U\r \Bbb R$ is a positive function $h=h(x,y)$.  Let $\om_{\Sigma}=h^2dx\w dy$ be the volume form of $\Sigma=(U,g)$. Let $M=U\times N$, where $N=\{(z,t)\in \Bbb R^2:z<0\}$. Define a metric $g$  on $M$  by $g(X,Y)=(\sinh a z)^2g_{\Sigma}(X,Y)+\th_3(X)\th_3(Y)+\th_4(X)\th_4(Y)$,   where $$\gather \th_3=\sinh 2azdt-(-\sin 2at\cosh 2az H(x,y)+\sinh 2az l_2(x,y))dx  \\-(\cos 2at \cosh 2az H(x,y)+\sinh 2azn_2(x,y))dy\endgather$$and
  $$\th_4=dz-\cos 2atH(x,y)dx-\sin 2atH(x,y)dy$$ and the function $H$ satisfies the equation  $\Delta \ln H=(\ln H)_{xx}+(\ln H)_{yy}=2a^2h^2+4a^2H^2$ on $U$, $l_2=\frac1{2a}(\ln H)_y, n_2=-\frac1{2a}(\ln H)_x$.  Then   $(M,g)$ admits a K\"ahler structure $\bJ$ with the K\"ahler form $$\bO = (\sinh az)^2\om_{\Sigma} +\th_4\w \th_3$$  and a Hermitian structure $J$ with the K\"ahler form $\Om=(\sinh az)^2\om_{\Sigma}+ \th_3\w \th_4$.  The Ricci tensor of $(M,g)$ is $J$-invariant and $J$ is not locally conformally K\"ahler.  The Lee form of $(M,g,J)$ is $\th=-\a\th_4$,  where $\a=-2a\coth az$.  }

\medskip

{\it Proof.} Take a coordinate system such that $E_1=\frac1f\p_x+k\p_z+l\p_t,E_2=\frac1f\p_y+m\p_z+n\p_t,E_3=\frac1\beta\p_t,E_4=\p_z$. Then  $\th_1=fdx, \th_2=fdy, \th_4=dz-(fk)dx-(fm)dy,\th_3=\beta dt-(\beta lf)dx-(\beta nf)dy$.  Let  $\a=-2a\coth az,\beta=\sinh 2az$. Then

$$E_1\ln\a=-\frac {2ak}{\sinh 2az}, \  \   E_2\ln\a=-\frac {2am}{\sinh 2az}.\tag5.5$$
We have $[E_1,E_4]=\frac{f_z}{f^2}\p_x-k_z\p_z-l_z\p_t=-\frac\a{2f}\p_x-\frac\a2k\p_z-\frac\a2 l\p_t-\frac {2am}{\sinh^22az}\p_t$. Hence $k_z=-a\coth az k,
l_z=-a\coth azl+\frac{2am}{\sinh^22az}, f_z=a\coth azf$.  This implies   $f=\sinh azh(x,y), k=\frac{k_1(x,y,t)}{\sinh az}$.  On the other hand, $[E_1,E_3]=-\frac1\beta( k_t\p_z+ l_t\p_t)-k\frac{\beta_z}{\beta^2}\p_t=\frac {2am}{\sinh 2az}\p_z$  and
$$l_t= -2ak\coth 2az, k_t=-2am.\tag 5.6$$
This yields
$$l=-m\coth 2az+\frac{l_1(x,y)}{\sinh a z}.\tag5.7$$

Similarly $[E_2,E_3]=-m\frac{\beta_z}{\beta^2}\p_t-\frac1{\beta} (m_t\p_z+ n_t\p_t)=-\frac {2ak}{\sinh 2az}\p_z$, and consequently,
$$m_t=2ak,n_t= -2am\coth 2az.$$
Since $[E_2,E_4]=-m_z\p_z-n_z\p_t+\frac{f_z}{f^2}\p_y=-\frac\a{2f}\p_y-\frac12\a m\p_z-\frac12\a n\p_t+\frac{2a k}{\sinh^22az}\p_t$, we get
$m_z=-a\coth az m$.  Hence  $m=\frac{m_1(x,y,t)}{\sinh az}$ and
$$n=k\coth 2az+\frac{n_1(x,y)}{\sinh az}.\tag5.8$$
Let us take  $kf=\cos 2atH(x,y), mf=\sin 2at H(x,y)$,  where  $f=\sinh azh(x,y)$ and $\th_1\w\th_2=(\sinh az)^2h^2dx\w dy$.  Then
$$\gather lf\beta=-\sin 2at\cosh 2azH(x,y)+\sinh 2azl_2(x,y),\\nf\beta=\cosh 2az\cos 2atH(x,y)+\sinh 2azn_2(x,y).\endgather$$  Hence
$$\gather \th_3=\sinh 2azdt-(-\sin 2at\cosh 2azH(x,y)+\sinh 2azl_2(x,y))dx\\-(\cos 2at\cosh 2azH(x,y)+\sinh 2azn_2(x,y))dy\endgather$$
and
$$\th_4=dz-\cos 2atH(x,y)dx-\sin 2atH(x,y)dy.$$

Now we prove that
$d\th_3=-\a\th_1\w\th_2-E_2\ln\a\th_1\w\th_4+E_1\ln\a\th_2\w\th_4+(-E_4\ln\a+\a)\th_3\w\th_4$ and
$d\th_4=E_2\ln\a\th_1\w\th_3-E_1\ln\a\th_2\w\th_3$  if $l_2=\frac1{2a}(\ln H)_y, n_2=-\frac1{2a}(\ln H)_x$, where  $\Delta\ln H=4a^2H^2+2a^2h^2$.

In fact,  $$\gather-\a\th_1\w\th_2-E_2\ln\a\th_1\w\th_4+E_1\ln\a\th_2\w\th_4+(-E_4\ln\a+\a)\th_3\w\th_4\\=2a\coth az f^2dx\w dy +\frac{2amf}{\sinh 2az}dx\w(dz-\sin 2atHdy)-\frac{2akf}{\sinh 2az}dy\\\w(dz-\cos 2atHdx)-\frac{2a\cosh 2az}{\sinh 2az}(\sinh 2azdt\\ +\cosh 2az \sin 2atH dx-\sinh 2azl_2 dx-\cos  2at\cosh 2az H dy -\sinh 2azn_2 dy)\\ \w(dz-\cos 2atHdx-\sin 2atHdy)=a\sinh 2azh^2dx\w dy\\-\frac{2a}{\sinh 2az}\sin^2 2atH^2dx\w dy +\frac{2a}{\sinh 2az}\cos^22atH^2dy\w dx\\ -\frac{2a\cos 2atH}{\sinh 2az}dy\w dz+\frac{2a\sin 2atH}{\sinh 2az}dx\w dz-\frac {2a\cosh 2az}{\sinh 2az}(\sinh 2azdt\w dz\\-\sinh 2az\cos 2atHdt\w dx-\sinh 2az\sin 2atHdt\w dy\\+\cosh 2az\sin 2atHdx\w dz-\sinh 2azl_2dx\w dz-\cosh 2az\sin^22atH^2dx\w dy\\ +\sinh 2azl_2\sin 2atHdx\w dy-\cosh 2az\cos 2atHdy\w dz\\-\sinh 2azn_2dy\w dz+\cosh 2az\cos^22atH^2dy\w dx +n_2\sinh 2az\cos 2atHdy\w dx). \endgather$$

On the other hand,
$$\gather d\th_3= 2a\cosh 2azdz\w dt +2a\sinh 2az\sin 2atHdz\w dx+\\2a\cosh 2az\cos 2atHdt\w dx+ \cosh 2az\sin 2atH_ydy\w dx\\-2al_2\cosh 2azdz\w dx-\sinh 2azl_{2y}dy\w dx +2a\cosh 2az\sin2a t Hdt\w dy\\-2a\cos 2at\sinh 2azHdz\w dy-\cos 2at\cosh 2azH_xdx\w dy-2a\cosh 2azn_2dz\w dy\\-\sinh 2azn_{2x}dx\w dy.\endgather$$

It is clear that $\th_3$ satisfies (2.3) if $l_2=\frac1{2a}(\ln H)_y, n_2=-\frac1{2a}(\ln H)_x$ and  $\Delta\ln H=2a^2h^2+4a^2H^2$, where $\Delta f=f_{xx}+f_{yy}$.

Now  $$\gather E_2\ln\a\th_1\w\th_3-E_1\ln\a\th_2\w \th_3\\=-2a\frac{\sin 2at H}{\sinh 2az}dx\w(\sinh 2az dt-\cos 2at\cosh 2az Hdy-\sinh 2azn_2dy)\\+\frac{2a\cos 2at H}{\sinh 2az} dy\w(\sinh 2azdt+\cosh 2az\sin 2at Hdx-\sinh 2azl_2dx)\\=2a\big(-\sin 2atHdx\w dt+\sin 2atHn_2dx\w dy +\frac{\cos 2at\sin 2at\cosh 2azH^2}{\sinh 2az}dx\w dy\\+\cos 2at Hdy\w dt+\frac{\sin 2at \cos 2at\cosh 2azH^2}{\sinh 2az}dy\w dx-\cos 2atHl_2dy\w dx\big)\\=-2a\sin 2atHdx\w dt+2a\sin 2atHn_2dx\w dy+2a\cos 2at Hdy\w dt\\-2a\cos 2atHl_2dy\w dx.\endgather$$

On the other hand,

$d\th_4=2a\sin 2atHdt\w dx-\cos 2atH_ydy\w dx-2a\cos 2atHdt\w dy-\sin 2atH_xdx\w dy$.

It is clear that $\th_4$ satisfies (2.3) if $l_2=\frac1{2a}(\ln H)_y, n_2=-\frac1{2a}(\ln H)_x$.

Note that   $d(\th_1\w\th_2)=d((\sinh az)^2h^2dx\w dy)=2a\sinh az\cosh azh^2dz\w dx\w dy=-(-2a\coth az)dz\w\th_1\w\th_2=-\a\th_4\w\th_1\w\th_2$. From (2.3) it is also clear that  $d(\th_3\w\th_4)=-\a\th_4\w\th_1\w\th_2$.   Hence  $d\bO=0$, where $\bO(X,Y)=g(\bJ X,Y)$ and $\bO=\th_1\w\th_2-\th_3\w\th_4$. Hence in view of Lemma F, $(M,g,\bJ)$ is a K\"ahler surface and the Lee form of $(M,g,J)$ is $\th=-\a\th_4$.  As above one can show  that  $(M,g,\bJ)$  is a QCH K\"ahler surface.$\k$

\bigskip
{\bf Theorem  4.}  {\it  Let  $U\subset \Bbb R^2$ and let $g_{\Sigma}=h^2(dx^2+dy^2)$ be a Riemannian metric on $U$, where $h:U\r \Bbb R$ is a positive function $h=h(x,y)$.  Let $\om_{\Sigma}=h^2dx\w dy$ be the volume form of $\Sigma=(U,g)$. Let $M=U\times N$, where $N=\{(z,t)\in \Bbb R^2:z<0\}$. Define a metric $g$  on $M$  by $g(X,Y)=(\cosh a z)^2g_{\Sigma}(X,Y)+\th_3(X)\th_3(Y)+\th_4(X)\th_4(Y)$,   where $$\gather\th_3=\sinh 2azdt-(\cos 2at\cosh 2azH(x,y)+\sinh 2azl_2(x,y))dx  \\-(-\cosh 2az\sin 2atH(x,y)+\sinh 2azn_2(x,y))dy\endgather$$and
  $$\th_4=dz-\sin 2atH(x,y)dx-\cos 2atH(x,y)dy$$ and the function $H$ satisfies the equation  $\Delta \ln H=(\ln H)_{xx}+(\ln H)_{yy}=-2a^2h^2+4a^2H^2$ on $U$, $n_2=\frac1{2a}(\ln H)_x, l_2=-\frac1{2a}(\ln H)_y$.  Then   $(M,g)$ admits a K\"ahler structure $\bJ$ with the K\"ahler form $$\bO = (\cosh az)^2\om_{\Sigma} +\th_4\w \th_3$$  and a Hermitian structure $J$ with the K\"ahler form $\Om=(\cosh az)^2\om_{\Sigma}+ \th_3\w \th_4$.  The Ricci tensor of $(M,g)$ is $J$-invariant and $J$ is not locally conformally K\"ahler.  The Lee form of $(M,g,J)$ is $\th=-\a\th_4$,  where   $\a=-2a\tanh az$ .  }

\medskip

{\it Proof.} Take a coordinate system such that $E_1=\frac1f\p_x+k\p_z+l\p_t,E_2=\frac1f\p_y+m\p_z+n\p_t,E_3=\frac1\beta\p_t,E_4=\p_z$. Then  $\th_1=fdx, \th_2=fdy, \th_4=dz-(fk)dx-(fm)dy,\th_3=\beta dt-(\beta lf)dx-(\beta nf)dy$.  Let  $\a=-2a\tanh az,\beta=\sinh 2az$. Then

$$E_1\ln\a=\frac {2ak}{\sinh 2az}, \  \   E_2\ln\a=\frac {2am}{\sinh 2az}.\tag5.8$$
We have $[E_1,E_4]=\frac{f_z}{f^2}\p_x-k_z\p_z-l_z\p_t=-\frac\a{2f}\p_x-\frac\a2k\p_z-\frac\a2 l\p_t+\frac {2am}{\sinh^22az}\p_t$. Hence $k_z=-a\tanh az k,
l_z=-a\tanh azl-\frac{2am}{\sinh^22az}, f_z=a\tanh azf$.  This implies   $f=\cosh azh(x,y), k=\frac{k_1(x,y,t)}{\cosh az}$.  On the other hand, $[E_1,E_3]=-\frac1\beta( k_t\p_z+ l_t\p_t)-k\frac{\beta_z}{\beta^2}\p_t=-\frac {2am}{\sinh 2az}\p_z$  and
$$l_t= -2ak\coth 2az, k_t=2am.\tag 5.6$$
This yields
$$l=m\coth 2az+\frac{l_1(x,y)}{\cosh a z}.\tag5.9$$

Similarly $[E_2,E_3]=-m\frac{\beta_z}{\beta^2}\p_t-\frac1{\beta} (m_t\p_z+ n_t\p_t)=\frac {2ak}{\sinh 2az}\p_z$, and consequently,
$$m_t=-2ak,n_t= -2am\coth 2az.$$
Since $[E_2,E_4]=-m_z\p_z-n_z\p_t+\frac{f_z}{f^2}\p_y=-\frac\a{2f}\p_y-\frac12\a m\p_z-\frac12\a n\p_t-\frac{2a k}{\sinh^22az}\p_t$, we get
$m_z=-a\tanh az m$.  Hence  $m=\frac{m_1(x,y,t)}{\cosh az}$ and
$$n=-k\coth 2az+\frac{n_1(x,y)}{\cosh az}.\tag5.8$$
Let us take  $kf=\sin 2atH(x,y), mf=\cos 2at H(x,y)$,  where  $f=\cosh azh(x,y)$ and $\th_1\w\th_2=(\cosh az)^2h^2dx\w dy$.  Then
$$\gather lf\beta=\cos 2at\cosh 2azH(x,y)+\sinh 2azl_2(x,y),\\nf\beta=-\cosh 2az\sin 2atH(x,y)+\sinh 2azn_2(x,y).\endgather$$  Hence
$$\gather \th_3=\sinh 2azdt-(\cos 2at\cosh 2azH(x,y)+\sinh 2azl_2(x,y))dx  \\-(-\cosh 2az\sin 2atH(x,y)+\sinh 2azn_2(x,y))dy\endgather$$
and
$$\th_4=dz-\sin 2atH(x,y)dx-\cos 2atH(x,y)dy.$$

Now we prove that
$d\th_3=-\a\th_1\w\th_2-E_2\ln\a\th_1\w\th_4+E_1\ln\a\th_2\w\th_4+(-E_4\ln\a+\a)\th_3\w\th_4$ and
$d\th_4=E_2\ln\a\th_1\w\th_3-E_1\ln\a\th_2\w\th_3$  if $l_2=-\frac1{2a}(\ln H)_y, n_2=\frac1{2a}(\ln H)_x$, where  $\Delta\ln H=4a^2H^2-2a^2h^2$.

In fact,  $$\gather-\a\th_1\w\th_2-E_2\ln\a\th_1\w\th_4+E_1\ln\a\th_2\w\th_4+(-E_4\ln\a+\a)\th_3\w\th_4\\=2a\tanh az f^2dx\w dy -\frac{2amf}{\sinh 2az}dx\w(dz-\cos 2atHdy)+\frac{2akf}{\sinh 2az}dy\\\w(dz-\sin 2atHdx)-\frac{2a\cosh 2az}{\sinh 2az}(\sinh 2azdt\\ -\cosh 2az \cos 2atH dx-\sinh 2azl_2 dx+\sin  2at\cosh 2az H dy -\sinh 2azn_2 dy)\\ \w(dz-\sin 2atHdx-\cos 2atHdy)=a\sinh 2azh^2dx\w dy\\+\frac{2a}{\sinh 2az}\cos^2 2atH^2dx\w dy -\frac{2a}{\sinh 2az}\sin^22atH^2dy\w dx\\ +\frac{2a\sin 2atH}{\sinh 2az}dy\w dz-\frac{2a\cos 2atH}{\sinh 2az}dx\w dz-\frac {2a\cosh 2az}{\sinh 2az}(\sinh 2azdt\w dz\\-\sinh 2az\sin 2atHdt\w dx-\sinh 2az\cos 2atHdt\w dy\\-\cosh 2az\cos 2atHdx\w dz-\sinh 2azl_2dx\w dz+\cosh 2az\cos^22atH^2dx\w dy\\ +\sinh 2azl_2\cos 2atHdx\w dy+\cosh 2az\sin 2atHdy\w dz\\-\sinh 2azn_2dy\w dz-\cosh 2az\sin^22atH^2dy\w dx +n_2\sinh 2az\sin 2atHdy\w dx). \endgather$$

On the other hand,
$$\gather d\th_3= 2a\cosh 2azdz\w dt -2a\sinh 2az\cos 2atHdz\w dx+\\2a\cosh 2az\sin 2atHdt\w dx- \cosh 2az\cos 2atH_ydy\w dx\\-2al_2\cosh 2azdz\w dx-\sinh 2azl_{2y}dy\w dx +2a\cosh 2az\cos2a t Hdt\w dy\\+2a\sin 2at\sinh 2azHdz\w dy+\sin 2at\cosh 2azH_xdx\w dy-2a\cosh 2azn_2dz\w dy\\-\sinh 2azn_{2x}dx\w dy.\endgather$$

It is clear that $\th_3$ satisfies (2.3) if $n_2=\frac1{2a}(\ln H)_x, l_2=-\frac1{2a}(\ln H)_y$ and  $\Delta\ln H=-2a^2h^2+4a^2H^2$, where $\Delta f=f_{xx}+f_{yy}$.

Note that  $$\gather -E_2\ln\a\th_1\w\th_3+E_1\ln\a\th_2\w \th_3\\=-2a\frac{\cos 2at H}{\sinh 2az}dx\w(\sinh 2azdt+\cosh 2az\sin 2at Hdy-\sinh 2azn_2dy)\\+\frac{2a\sin 2at H}{\sinh 2az} dy\w(\sinh 2az dt-\cos 2at\cosh 2az Hdx-\sinh 2azl_2dx)\\=2a\big(-\cos 2atHdx\w dt+\cos2atHn_2dx\w dy -\frac{\cos 2at\sin 2at\cosh 2azH^2}{\sinh 2az}dx\w dy\\+\sin 2at Hdy\w dt-\frac{\sin 2at \cos 2at\cosh 2azH^2}{\sinh 2az}dy\w dx-\sin 2atHl_2dy\w dx\big)\\=-2a\cos 2atHdx\w dt+2a\cos 2atHn_2dx\w dy+2a\sin 2at Hdy\w dt\\-2a\sin 2atHl_2dy\w dx.\endgather$$

On the other hand,

$d\th_4=-2a\cos 2atHdt\w dx-\sin 2atH_ydy\w dx+2a\sin 2atHdt\w dy-\cos 2atH_xdx\w dy$.

It is clear that $\th_4$ satisfies (2.3) if $n_2=\frac1{2a}(\ln H)_x, l_2=-\frac1{2a}(\ln H)_y$.

Now we have   $d(\th_1\w\th_2)=d((\cosh az)^2h^2dx\w dy)=2a\sinh az\cosh azh^2dz\w dx\w dy=-(-2a\tanh az)dz\w\th_1\w\th_2=-\a\th_4\w\th_1\w\th_2$. From (2.3) it is also clear that  $d(\th_3\w\th_4)=-\a\th_4\w\th_1\w\th_2$.   Hence  $d\bO=0$, where $\bO(X,Y)=g(\bJ X,Y)$ and $\bO=\th_1\w\th_2-\th_3\w\th_4$. Hence in view of Lemma F, $(M,g,\bJ)$ is a K\"ahler surface and the Lee form of $(M,g,J)$ is $\th=-\a\th_4$.  As above one can show  that  $(M,g,\bJ)$  is a QCH K\"ahler surface.$\k$
\medskip
{\it Remark.}   Note that the generalized Calabi type K\"ahler surfaces which are not of Calabi type are  fiber bundles over $\Sigma$  and the fibers which are leafs of the foliation $\De$ have constant sectional curvature  :  $0$  in the semi-symmetric case,   positive $4a^2$ if $\a=2a\tan az$ and negative $-4a^2$ if $\a=-2a\coth az$ or $\a=-2a\tanh az$.

\par
\bigskip

\centerline{\bf References.}

\cite{A-C-G} V. Apostolov, D.M.J. Calderbank, P. Gauduchon {\it
The geometry of we\-akly self-dual K\"ahler surfaces},  Compos. Math.
135, 279-322, (2003)
\par
\medskip
\cite{A-G} V. Apostolov and P. Gauduchon, {\it The Riemannian Goldberg-Sachs theorem}, Int. J. Math. {\bf 8}, (1997), 421-439.
par
\medskip
\cite{D-M-1} A. Derdzi\'nski, G. Maschler {\it Special
K\"ahler-Ricci potentials on compact K\" ahler manifolds}, J.
reine angew. Math. 593 (2006), 73-116.
 \par
\medskip
\cite{D-M-2}  A. Derdzi\'nski, G.  Maschler {\it Local
classification of conformally-Einstein  K\"ahler metrics in higher
dimensions}, Proc. London Math. Soc. (3) 87 (2003), 779-819.
\par
\medskip
\cite{D-T}   Maciej Dunajski and Paul Tod {\it  Four-dimensional metrics conformal to K\"ahler }
Mathematical Proceedings of the Cambridge Philosophical Society, Volume 148, (2010), 485-503

\par
\medskip
\cite{J-1} W. Jelonek {\it K\"ahler surfaces with quasi constant holomorphic curvature}, Glasgow Math. J. 58, (2016), 503-512.
\par
\medskip
\cite{J-2} W. Jelonek {\it Semi-symmetric  K\"ahler surfaces }, Colloq. Math.  148, (2017), 1-12.
\par
\medskip
\cite{J-3} W. Jelonek {\it Complex foliations and  K\"ahler QCH surfaces },  Colloq. Math.  156, (2019), 229-242.
\par
\medskip
\cite{J-4} W. Jelonek {\it Einstein Hermitian and anti-Hermitian 4-manifolds }, Ann.  Po\-lon. Math. 81, (2003), 7-24.
\medskip
\cite{J-5} W. Jelonek {\it  QCH  K\"ahler surfaces II},  Journal of Geometry and Physics, (2020) 103735,
https://doi.org/10.1016/j.geomphys.2020.103735

\par
\medskip
Authors adress:

 Institute of Mathematics

 Cracow  University of Technology

 Warszawska 24

 31-155 Krak\'ow,  POLAND.

\enddocument